\documentclass[5p,times]{elsarticle}
\usepackage{amssymb}
\usepackage{amsthm}
\usepackage{makecell}
\usepackage{amsmath}
\usepackage{color,xcolor}
\newtheorem{remark}{Remark}
\usepackage[utf8]{inputenc}
\usepackage{booktabs, caption, makecell}

\usepackage{threeparttable}

\usepackage{multirow}
\usepackage[figuresright]{rotating}
\usepackage{amsmath}
\newcommand{\R}{\mathbb{R}}
\usepackage{tablefootnote}
\newcommand{\Z}{\mathbb{Z}}
\DeclareMathOperator*{\minimize}{minimize}
\DeclareMathOperator*{\maximize}{maximize}
\DeclareMathOperator*{\subt}{subject \ to}

\newtheorem{problem}{Problem}
\usepackage{amsmath, amssymb, amsfonts, amsthm}
\usepackage[linesnumbered, lined, boxed]{algorithm2e}

\begin{document}

\begin{frontmatter}

\title{Taking the human out of decomposition-based optimization via artificial intelligence: Part II. Learning to initialize}

\author{Ilias Mitrai}
\author{Prodromos Daoutidis\fnref{label2}}
\fntext[label2]{Corresponding author}
\ead{daout001@umn.edu}
\address{Department of Chemical Engineering and Materials Science, University of Minnesota, Minneapolis, MN 55455}

\begin{abstract} 
The repeated solution of large-scale optimization problems arises frequently in process systems engineering tasks. Decomposition-based solution methods have been widely used to reduce the corresponding computational time, yet their implementation has multiple steps that are difficult to configure. We propose a machine learning approach to learn the optimal initialization of such algorithms which minimizes the computational time. Active and supervised learning is used to learn a surrogate model that predicts the computational performance for a given initialization. We apply this approach to the initialization of Generalized Benders Decomposition for the solution of mixed integer model predictive control problems. The surrogate models are used to find the optimal number of initial cuts that should be added in the master problem. The results show that the proposed approach can lead to a significant reduction in solution time, and active learning can reduce the data required for learning.
\end{abstract}

\end{frontmatter}

\section{Introduction}
Decomposition-based optimization algorithms have been used widely to solve complex and large-scale optimization problems in a broad range of applications in chemical engineering, such as production scheduling and planning, supply chain management, mixed integer optimal control, and real-time operation. These algorithms exploit the underlying structure of a problem and decompose it into a number of easier-to-solve subproblems. Typical examples include Benders \cite{Benders1962} and Generalized Benders Decomposition \cite{geoffrion1972generalized}, Lagrangean decomposition \cite{guignard1987lagrangean}, Alternating Direction Method of Multipliers (ADMM) \cite{boyd2011distributed}, and cross decomposition \cite{van1983cross}. 

Despite the wide success of these algorithms, their efficiency over monolithic methods is not known a-priory and their implementation, especially in an online setting, is challenging due to the many steps involved and the underlying computational complexity of every step. Automatically determining when to use and how to implement a decomposition-based solution algorithm can simplify their implementation and reduce the solution time for complex optimization problems. In the companion paper \cite{mitrai2023whentodec} we proposed a graph classification approach to determine when to use a decomposition-based solution algorithm. In this paper, we focus on the implementation of decomposition-based solution methods.

The application of a decomposition-based solution algorithm has three steps: (1) problem decomposition, (2) coordination scheme, and (3) initialization of the algorithm. In the first step, the original problem is decomposed into a number of easier-to-solve subproblems. This step requires knowledge of the underlying structure of the problem. 
Automatic decomposition approaches either represent the problem as a graph and employ methods from network science to learn the underlying structure \cite{allman2019decode,mitrai2020decomposition,mitrai2022stochastic,mitrai2021iecr,jalving2018graph} or use machine learning (ML) to analyze the computational efficiency of different decompositions \cite{basso2020random, basso2023data}. 
The coordination step determines the exchange of information between the different subproblems. For distributed algorithms, the coordination determines the update of the dual variables whereas for hierarchical algorithms it involves the addition of cuts. 
Finally, the last step is the initialization of the algorithm. Unlike monolithic algorithms which may require an initial guess for the variables, decomposition-based algorithms require additional information regarding the values of the dual variables for distributed algorithms or an initial set of cuts for hierarchical algorithms. The configuration of these steps can have an important effect on the computational performance of a decomposition-based solution algorithm, however, determining the best configuration for a given problem is nontrivial. 

These steps can be viewed as hyperparameters of the solution algorithm. Therefore, given an optimization problem and a decomposition-based solution algorithm, one must find the optimal values of the algorithm parameters such that a desired performance function, e.g., the solution time, is optimized. Formally this is known as the algorithm configuration problem and is stated as follows \cite{eggensperger2019pitfalls,  schede2022survey}:
\begin{problem}
    \normalfont (Algorithm configuration) Given an optimization problem $P$, an optimization algorithm $\alpha$ with parameter $n \in \mathcal{N}$, and a performance function $m: \mathcal{P} \times \Pi \mapsto \mathcal{M}$, determine the optimal values of the parameters $n^*$ that optimize $m$ \\
    \begin{equation}
    n^* \in \arg \min_{n \in \mathcal{N}} \ \ m(n,P).
\end{equation}
\end{problem}
The algorithm configuration problem has three components, the optimization problem $P$ which belongs in some class of optimization problems $\mathcal{P}$, the parameter space, i.e., all possible values of the parameters represented by $\mathcal{N}$, and the performance space $\mathcal{M}$ which is a metric to compare the different configurations. We note that the version of the problem presented above is known as the per-instance algorithm configuration problem since the optimal parameters are determined only for a given problem $P$. Alternatively, one can find the optimal values of the parameters for a class or set of optimization problems. The algorithm configuration can be either static, i.e., the parameters of the algorithms remain fixed during the solution process, or dynamic, where the parameters adapt as the solution procedure evolves. The solution of the algorithm configuration problem is challenging since optimization solvers, either monolithic or decomposition-based, have multiple algorithmic steps and each step can have different parameters. Furthermore, optimization solvers employ a number of heuristics; although these can accelerate the solution on average, their efficacy for a given problem is not known a-priory. Therefore, finding the optimal configuration of an optimization solver is a challenging black-box optimization problem since the number of possible combinations of parameters can be very large, their optimal values may change significantly for different classes of optimization problems, and the evaluation of a given configuration can be computationally expensive.  

Automated algorithm configuration approaches aim to either solve the algorithm configuration problem in a computationally efficient way or approximate its solution. In the former approach,   Bayesian optimization and derivative-free methods have been employed to tune optimization algorithms by optimizing directly the black-box performance function \cite{liu2019tuning,chen2011random,hutter2009paramils,hutter2011sequential,hutter2010automated}. In the latter approach,  ML tools have been used to approximate the performance function $m$ with a surrogate one $\hat{m}$ and then identify the optimal values of the parameters $n^*$. These methods have been traditionally applied for tuning monolithic solvers, either by considering all the parameters simultaneously \cite{iommazzo2020learning} or specific algorithmic steps such as branching \cite{lodi2017learning, khalil2016learning, balcan2018learning, di2016dash, gupta2020hybrid, gupta2022lookback, liu2022learning}, cutting plane methods \cite{tang2020reinforcement,huang2022learning, paulus2022learning}, estimating active constraints \cite{misra2022learning, bertsimas2021voice,bertsimas2022online}, and improving primal heuristics \cite{ding2020accelerating}. 
ML tools have also been used to tune decomposition-based solution methods. For hierarchical decomposition-based methods, such as column generation and Benders Decomposition, ML is used to aid the cut selection process using classification techniques \cite{morabit2021machine, jia2021benders,lee2020accelerating}, whereas for distributed algorithms, such as ADMM, machine learning is used to determine the values of the dual variables and penalty parameters \cite{biagioni2020learning, zeng2022reinforcement}. In both cases, ML is used in the coordination step of these algorithms.

These ML approaches are based on a handcrafted feature representation of an optimization problem which is the input to the surrogate model $\hat{m}$ \cite{smith2012measuring, hutter2014algorithm,bengio2021machine,chen2021learning}, and involve two steps; the first is offline where multiple problems are solved for different values of the parameters in order to create a training dataset used to learn a surrogate model for the performance function. Once the surrogate model is trained, it is used online to identify the best set of parameters for a given problem. The main limitation of these ML approaches is data availability since generating a large training data set can be computationally expensive. 

In this work, we focus on the \textit{initialization} of decomposition-based solution methods. Specifically, we focus on cutting plane-based hierarchical solution methods, such as Benders and Generalized Benders Decomposition. These algorithms are based on the observation that if a subset of the variables, called complicating variables, is fixed then the problem is either easier to solve or has a special structure. Thus the original problem is decomposed into two problems; a master problem which considers the complicating variables and a subproblem which considers the non-complicating variables and whose solution depends on the values of the complicating variables. The solution of the master problem and the subproblem is coordinated via the addition of optimality and feasibility cuts, which inform the master problem about the bounds and the feasibility of the problem respectively. In the standard application of the algorithm, the master problem is initially solved without any cuts. Since the cuts contain information about the effect of the complicating variables on the subproblem, the addition of cuts in the first iteration can potentially lead to a reduction in the computational time since fewer iterations might be necessary.  However, the addition of a large number of cuts may increase the computational complexity of the master problem which in turn may increase the solution time. Hence it is important to identify the number of cuts that balances the amount of information added to the master problem with the increase in computational complexity. This balance is especially important for online applications where an optimization problem is solved repeatedly to compensate for updated process information. Depending on the application, either only the parameters of the optimization problem change, such as in model predictive control applications \cite{tang2018optimal, mitrai2023mipmpc}, or the parameters and the number of variables and constraints can change, such as in online scheduling applications \cite{mitrai2022adaptive,risbeck2019unification}. 

In this paper, we propose an ML approach to learn how to initialize Generalized Benders Decomposition. The proposed approach has two steps. In the first,  ML is used to learn a surrogate model that estimates the solution time of the optimization problem for a given number of cuts added to the master problem in the first iteration of the algorithm. In the second step, the surrogate model is used online to determine the optimal number of cuts that should be added in the master problem in the first iteration. 
We apply the proposed approach to a case study on the real-time operation of process systems, where a mixed integer economic model predictive control problem is solved. Specifically, we assume that the system is an isothermal continuously stirred tank reactor (CSTR) that can manufacture a number of products, and multiple disturbances can affect the operation of the system. The mixed integer economic model predictive control problem is solved using a hybrid multicut Generalized Benders Decomposition proposed in \cite{mitrai2022multicut}, initialized using the proposed approach. The results show that (1) the optimal initialization can be achieved in an automated way without human intervention, (2) the proper initialization can lead to a significant reduction in solution time, and (3) active learning can guide the learning process either during the initial development of such frameworks or for cases where generating the training dataset is computationally expensive. 

The rest of the paper is organized as follows: In Section \ref{GBD intro} we present the Generalized Benders Decomposition algorithm and the different acceleration techniques, in Section \ref{alg conf prob} we pose the initialization of Generalized Benders Decomposition as an algorithm configuration problem, in Section~\ref{L2InitViaAL} we present the proposed approach, and in Section~\ref{application} we present the case study and the numerical results. Finally, in Section~\ref{concls}, given the results in the first part of this two-series paper \cite{mitrai2023whentodec} we present a unified framework for automated decomposition-based solution algorithm selection and configuration.

\section{Generalized Benders Decomposition} \label{GBD intro}
\subsection{Standard implemetation}
 We will assume that the following problem (denoted as $P$) must be solved:
\begin{equation} \label{full model}
    \begin{aligned}
    P (p):= \minimize_{z,x,y} \ \ & f_1(z,x;p_m) + f_2(x,y;p_s) \\
    \subt \ \ & g_1(z,x;p_m) \leq 0\\
              & h_1(z,x;p_m) = 0\\
              & g_2(x,y;p_s) \leq 0\\
              & h_2(x,y;p_s) =0\\
              & z \in \Z^{n_z}, x \in \R^{n_x^c}\times \Z^{n_x^d}, y \in \R^{n_y},
    \end{aligned}   
\end{equation}
where $p=[p_m, p_s]^\top$ are the parameters of the problem, and  $z \in \Z^{n_z}$, $x \in \R^{n_x^c}\times \Z^{n_x^d}$,$ y \in \R^{n_y}$ are the variables. The solution of this problem depends on the values of the parameters $p$. In this problem, we observe that if the variables $z,x$ are fixed, then the resulting problem is a continuous optimization problem that depends on the values of the variables $x$, which are the complicating variables, and the parameters $p_s$. Given this structure, we can apply Generalized Benders Decomposition, by assigning the $z,x$ variables and the associated constraints in the master problem and the other variables and constraints in the subproblem. Under this decomposition, parameters $p_m$ affect only the master problem and parameters $p_s$ affect only the subproblem. The subproblem is
\begin{equation} \label{sub alg 1}
    \begin{aligned}
    \mathcal{S} (x,p_s):= \minimize_{\bar{x},y} \ \ & f_2(\bar{x},y;p_s) \\
    \subt \ \ & g_2(\bar{x},y;p_s) \leq 0\\
              & h_2(\bar{x},y;p_s) =0\\
              & \bar{x} =x \ \ : \ \  \lambda\\
              & \bar{x} \in \R^{n_x^c+n_x^d}, y \in \R^{n_y},
    \end{aligned}   
\end{equation}
where $\lambda$ are the Lagrangean multipliers for the equality constraint $\bar{x}=x$. The solution of this problem depends on the values of the complicating variables $x$ and the parameters $p_s$, and for a given value of $x=\bar{x}$ the value function of the subproblem can be approximated as follows:
\begin{equation} \label{benders cut}
    \mathcal{S}(x, p_s) \geq \mathcal{S}(\bar{x}, p_s) - \bar{\lambda} (x - \bar{x}),
\end{equation}
where $\bar{\lambda}$ is equal to the value of the Lagrangean multiplier at the optimal solution of the subproblem when solved for $x = \bar{x}$. Note that we assume that the subproblem is always feasible for all values of $x$. The master problem is:
\begin{equation} \label{master}
    \begin{aligned}
    \mathcal{M}(\cdot):= \minimize_{z,x,\eta } \ \ & f_1(z,x;p_m) + \eta \\
    \subt \ \ & g_1(z,x;p_m) \leq 0\\
              & h_1(z,x;p_m) = 0\\
              & \eta \geq \mathcal{S}(\bar{x}^l,p_s) - \lambda^l (x-\bar{x}^l) \ \forall l \in \mathcal{L} \\
              & z \in \Z^{n_z}, x \in \R^{n_x^c}\times \Z^{n_x^d},
    \end{aligned}   
\end{equation}
where $\mathcal{M}(\cdot) = \mathcal{M}(p_m, \mathcal{L})$, $l$ is the iteration number and the set $\mathcal{L}$ denotes the index of the Benders cuts. The steps of GBD are presented in Algorithm 1. 

\begin{algorithm}
\caption{Generalized Benders Decomposition}\label{alg:benders algorithm}
\KwData{Optimization problem}
\KwResult{Upper, lower bound and variable values}
Set $UB = \infty$, $LB = - \infty$\;
Set tolerance and optimality gap (tol)\;
Initialize the algorithm\;
\While{$(UB-LB)/LB \geq \text{tol}/100$}{
Solve the master problem (Eq.~\ref{master}) and obtain $LB,x$\;
Solve the subproblem (Eq.~\ref{sub alg 1}) and obtain $y$\;
Add Benders cut (Eq.~\ref{benders cut})\;
Update the upper bound $f_1 (z,x;p_m)+f_2(x,y;p_s)$\;
}
\end{algorithm}

\subsection{Acceleration techniques for Benders decomposition}
The algorithm alternates between the solution of the master problem and the subproblem, therefore the computational performance depends on the complexity of the master problem and subproblem, i.e., the solution time per iteration, the number of infeasible subproblems that must be solved, and the quality of cuts that are generated during the solution. Two approaches can be followed to handle these issues. The first one is based on the theoretical aspects of the algorithm and the underlying geometry of the problem. Common strategies in this approach are problem reformulation and decomposition \cite{crainic2016partial, magnanti1981accelerating}, the addition of valid inequalities in the master problem to reduce the number of infeasible subprolems \cite{saharidis2011initialization}, multicut implementation \cite{you2013multicut} for stochastic optimization problems, cut generation and management \cite{magnanti1981accelerating, su2015computational, saharidis2010improving, pacqueau2012fast, varelmann2022decoupling}, and regularization/stabilization of the master problem \cite{ruszczynski1997accelerating, linderoth2003decomposition}. 
The second approach involves the use of ML techniques. For example, ML has been used to develop a classifier to determine which cuts should be added in the master problem during the multicut implementation of the algorithm for the solution of two-stage mixed integer stochastic optimization problems \cite{jia2021benders, lee2020accelerating}. ML has also been used to approximate the solution of the subproblem reducing the solution time for cases where the subproblem is computationally complex, such as two-stage stochastic optimization problems \cite{larsen2023fast} and mixed integer model predictive control problems \cite{mitrai2023mipmpc}. These acceleration methods, both the ones based on the underlying geometry and the ones using ML, reduce the computational time either by reducing the solution time per iteration or by reducing the number of iterations.

\section{Initialization of GBD as an algorithm configuration problem} \label{alg conf prob}

We will assume that problem $P(p)$ must be solved repeatedly given new values of the parameters $p$. The problem that we will address is the following: 
\begin{problem}
    \normalfont Given an optimization problem $P(p)$ determine the optimal number of cuts to add in the master problem in the first iteration of Generalized Benders Decomposition, such that the CPU time is minimized.
\end{problem}
We can pose this as a per-instance algorithm configuration problem as follows
\begin{equation}
    \minimize_{n \in \mathcal{N}_{cuts}} \ m( n , P(p)),
\end{equation}
where the optimization problem is $P(p)$, the parameter space $\mathcal{N}_{cuts}$ represents the cuts that can be added, and the performance function is the solution time $\mathcal{M}= \mathbb{R}_{+}$. The solution of the algorithm configuration problem for Generalized Benders Decomposition has three main challenges. The first, which is common in algorithm configuration problems, is that the performance function $m$ is not known a-priory. The second issue is related to the parameter space, i.e., all the cuts that can potentially be used to initialize the master problem. In general, selecting which cuts to use is a challenging problem that arises in the solution of many classes of problems such as quadratic \cite{baltean2019scoring, marousi2022acceleration}, and mixed integer linear programming problems \cite{morabit2021machine,chi2022deep}. Regarding Generalized Benders Decomposition, the number of cuts that can be evaluated depends on the number and type of complicating variables. For cases where the complicating variables are integer and $n_c$ cuts must be added, the cut selection process leads to a combinatorial optimization problem since one must select $n_c$ cuts from all possible cuts. The situation is more complex when the complicating variables are continuous, since in this case an infinite number of cuts can be added even for one complicated variable. Finally, all the parameters of the problem can change simultaneously, thus the cuts must be evaluated continuously as the parameter of the problem change. Hence, although the addition of cuts in the master problem might reduce the number of iterations required for convergence, the solution time might increase since multiple subproblems must be solved to evaluate the initial set of cuts.

\par
\begin{figure*}
    \centering
    \includegraphics[scale=0.18]{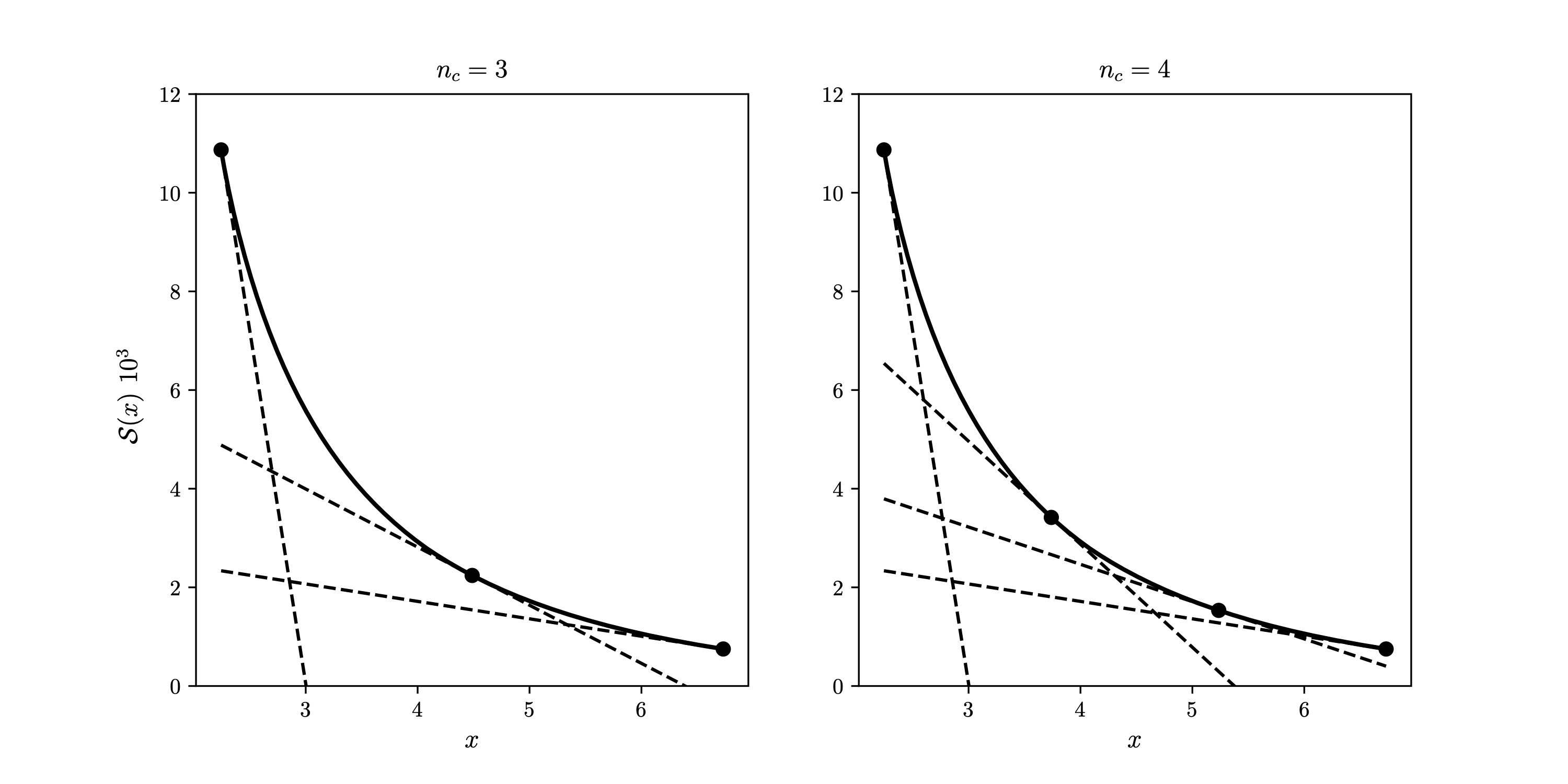}
    \caption{Domain discretization for the case study considered in Section~\ref{application} for a transition from product 1 to 2 for three and four number of cuts $(n_c)$. The solid line is the value function, $x\in [2.24, 6.73]$ is the complicating variable, the dotted lines are the value function approximations, i.e., Benders cuts, evaluated ad the points indicated by the dots. }
    \label{fig:cut discr scheme}
\end{figure*}

 We will assume that (1) the parameters of the subproblem do not change, (2) all the complicating variables are continuous, and (3) $n_{c}$ cuts can be added by discretizing the domain of the complicating variables ($x \in [x^{lb}, x^{ub}]$) into $n_{c}$ uniform points. The first assumption guarantees that the Benders cuts must be evaluated only once since even if the parameters of the master problem change, the Benders cuts are still valid estimators of the value function of the subproblem, since the parameters of the subproblem do not change. Therefore, the process of adding these cuts to the master problem does not affect the total solution time. The second and third assumptions determine the cut selection strategy. In this work, the number of cuts determines the number of points that will be used to approximate the value function of the subproblem via Benders cuts. An example is presented in Fig.~\ref{fig:cut discr scheme}, where the value function corresponds to the system discussed in Section~\ref{application} and the approximation with three and four cuts is presented. This setting enables us to simplify the cut selection process for the case of continuous complicating variables. Notice that in this setting the cuts are uniformly distributed in the domain of the complicating variables. 

\section{Learning to initialize via supervised and active learning} \label{L2InitViaAL}
The goal is to learn a surrogate model $\hat{m}$ which approximates the solution time and can be used to identify the optimal number of cuts to add in the master problem, 
\begin{equation}
    \begin{aligned}
    n^* \in \arg \min_{n \in \mathbb{Z}_{+}} \ \hat{m}(n,\nu( \mathcal{P}(p))),
    \end{aligned}
\end{equation}
where $\nu(P(\textbf{p}))$ are some features of problem $P(p)$. We propose two ML approaches for learning the surrogate model. 

\subsection{Supervised learning approach}
The estimation of the parameters of a surrogate model requires data that capture the relation between the features of an optimization problem $\nu(P(p))$ and the number of cuts $n$ with the solution time. Such data can be generated using the procedure presented in Algorithm~\ref{alg:supervised learning}, where first random values of the parameters of the master problem $(p_m)$ are generated based on some underlying probability distribution. Next, the optimization problem $P(p_i)$ is solved for a fixed number of cuts $n_i$, and the solution time $y_i$ is obtained. Finally, the features of the problem $\nu(P(p_i))$ are obtained and the tuplet $(n_i, \nu(P(p_i)))$ and the solution time $y_i$ are stored. Once this procedure is completed, we obtain the dataset $ \mathcal{D} = \{ (n_i, \nu(P(p_i))), y_i\}_{i=1}^{N_{data}}$. This dataset can be used to learn the parameters of the surrogate model by optimizing some loss function such as the squared error between the model prediction and the data. Once the learning step is completed, the surrogate model can be used to learn how many cuts to add for online applications as presented in Algorithm~\ref{alg:online initialization}, where given new values of the parameters of the optimization problem, the features of the problem are obtained and the optimal number of cuts is computed by optimizing the surrogate model. Finally, the cuts are added to the master problem and Generalized Benders Decomposition is implemented.  

\begin{algorithm}[t]
\caption{Learning the relation between number of cuts and CPU time for a general optimization problem for continuous complicating variables} \label{alg:supervised learning}
\KwData{Optimization problem, number of data points $N_{data}$, number of discretization points $N_{cuts}$, upper and lower bounds for complicating variables $x \in [x^{lb}, x^{ub}]$, $\check{p} = [p_4,p_5,p_6]$}
\KwResult{Surrogate model $\hat{m}$}
$i=1$\;
\While{$i \leq N_{data}$}{
Generate parameters $p_1,p_2,p_3 \rightarrow p_i = [p_1,p_2,p_3, \check{p}]$\;
\For{$j=2:N_{cuts}$}{
Solve problem $\mathcal{P}(p_i)$ using $j$ cuts for each $x$\;
Obtain CPU time $y_i$\;
Obtain features of the problem $s_j = (\nu( \mathcal{P}(p_i)),j)$\;
Append data $\{s_j,y_j\}$\;
}
$i = i+1$\;
}
Using data $\{s_i,y_i\}_{i=1}^{N_{data}(N_{cuts}-1)}$ learn parameters of a surrogate model $\hat{m}$\;
\end{algorithm}

\subsection{Active learning approach}
The main limitation of the supervised learning approach is the computational time required to generate the training data, since for every value of the parameters $p_i$ and number of cuts $n_i$ the problem must be solved to obtain the solution time. This approach can be computationally expensive, even intractable for complex optimization problems. To resolve this we propose the application of active learning \cite{settles2009active}, a commonly used approach in ML tasks where the features of the data are known but obtaining their label is costly or time-consuming. Unlike supervised learning where all the data are available for training, in active learning the model itself determines which data should be labeled and thus be used for training the surrogate model. 

The active learning paradigm has three components. The first is the unlabeled data which can  be generated de novo (known as membership query synthesis) \cite{angluin1988queries}, can become available in an online setting (stream-based selective sampling) \cite{cohn1994improving}, or can be gathered at once (pool-based sampling) \cite{lewis1995sequential}. The second aspect is the query strategy, which determines which data should be labeled. This decision is taken by considering the informativeness of the available unlabeled data. Typical query strategies are uncertainty sampling \cite{lewis1995sequential}, query by committee \cite{seung1992query}, and expected model change \cite{settles2007multiple}. We refer the reader to \cite{settles2009active} for a detailed discussion of the different sampling methods. The last component is an oracle which generates the label for a given input. Typical example of an oracle is a human expert, a computer simulation or the outcome of an experiment. 

\begin{figure}
	\centering
	\includegraphics[trim = 0 00 0 0, clip,scale=0.35]{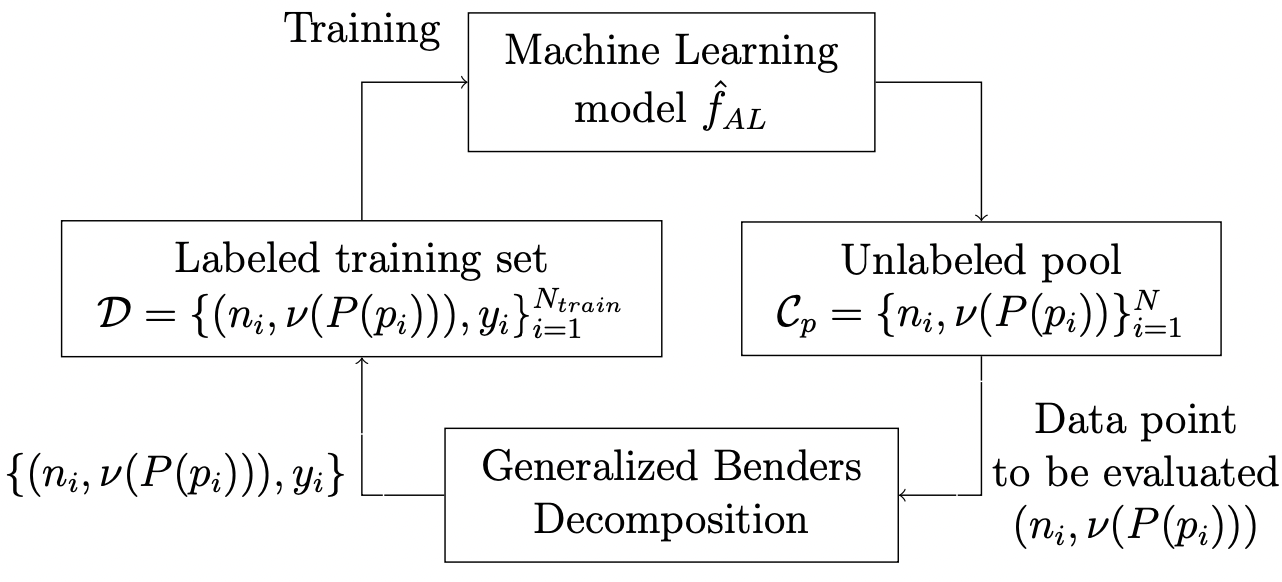}
	\caption{Learning to initialize Generalized Benders Decomposition via active learning framework}
    \label{fig: L2InitAL}
\end{figure}

In this work, we will use the pool-based active learning paradigm with uncertainty-based sampling, where the data point for which the model is the least certain is labeled. The basic steps of the application of active learning for learning the solution time of Generalized Benders Decomposition are presented in Fig.~\ref{fig: L2InitAL}.

\subsubsection{Generation of pool of labels and initial training set}
For the application of active learning, first, the pool of labels is created. We generate random values of the parameter $p_i$ and for every parameter the features of the optimization problem $\nu(P(p_i))$ are obtained and a number of cuts $n_i$ is selected. This forms the pool of features $\mathcal{C}_p = \{ s_i \}_{i=1}^{N_{pool}}$ $(s_i = \{ n_i, \nu(P(p_i))\})$. Next, a small set of training data is obtained by sampling $N_{initial}$ data points from the pool and evaluating the solution time $y_i$. This is the initial training set $\mathcal{C} = \{s_i,y_i\}_{i=1}^{N_{initial}}$ which we will refer to as labeled training set. Given the unlabeled pool $\mathcal{C}_p$ and the labeled dataset $\mathcal{C}$, the surrogate model considers all the data in the pool, and the data point $s = \{ n_s, \nu(P(p_s)) \}$ for which the prediction is the least certain about is selected for labeling. 

\subsubsection{Uncertainty-based sampling using Gaussian Processes}
The selection of the data point requires quantification of the prediction uncertainty. Gaussian Process Regression (GPR) is a non-parametric Bayesian approach that can provide uncertainty measures \cite{williams2006gaussian}. GPR is based on the Gaussian Process which is a stochastic process that defines a distribution over functions. Specifically, given observations of the input variables $X = [ x_1,...,x_N]$ and measurements of the output $Y = [ y_1,...,y_N]$, the relationship between $X$ and $Y$ is modelled as a Gaussian multivariate distribution. GPR seeks to learn a mapping $f: X \mapsto Y$, i.e., $y=f(x)$, with mean $m(x)$ and covariance $k(x,x')$ where ($\mathbb{E}$ is the expected value)
\begin{equation}
\begin{aligned}
    m(x) & = \mathbb{E}[ f(x) ] \\
    k(x,x') & = \mathbb{E} [ (f(x) - m(x)) (f(x')-m(x'))].
\end{aligned}   
\end{equation}
This is done under the assumption that the data are independent and the probability to observe an output given the observations can be factored over cases in the training set. The Gaussian Process, is written as $f(x) \sim \mathcal{GP}(m(x), k(x,x'))$. Different kernel functions $k(\cdot,\cdot)$ can be used, however in this paper we use the Matern kernel given by the following equation
\begin{equation}
    k(\cdot,\cdot) = \frac{1}{\Gamma(\nu) 2^{\nu-1}} \bigg(\frac{\sqrt{2\nu}}{\ell} d(\cdot,\cdot) \bigg)^{\nu} K_v \bigg( \frac{\sqrt{2 \nu}}{\ell} d(\cdot,\cdot) \bigg),
\end{equation}
where $d(\cdot,\cdot)$ is the Eucledian distance between features $s_i$ and $s_j$, $\Gamma(\cdot)$ is the gamma function, $K_v(\cdot)$ is the modified Bessel function, and $\ell, \nu$ are tunable hyperparameters. During training the parameters of the mean and kernel function are estimated based on the available data. This step estimates the posterior distribution over functions that best explain the data. This posterior distribution is also Gaussian and is used to make predictions for a new data point. We refer the reader to \cite{williams2006gaussian} for detailed explanation of these steps.

\subsubsection{Oracle}
The oracle is the Generalized Benders Decomposition algorithm which given the parameters of an optimization problem and a number of cuts, solves the optimization problem and returns the solution time. We note that although here we consider the standard version of Generalized Benders Decomposition as the oracle, in principle, other versions can be incorporated, such as multicut implementation, partial Benders decomposition etc. 

\subsubsection{Active learning loop}
The active learning loop is presented in Fig.~\ref{L2InitViaAL} and Algorithm~\ref{alg:AL algo}. The inputs are the pool of unlabeled data $\mathcal{C}_p$, the initial training set $C$, and the surrogate model $\hat{m}$ which is a Gaussian Process with Matern kernel. First, the model is trained using the initial training dataset. Next, the model is used to predict the solution time and uncertainty around the prediction for all the unlabeled data and identify the datapoint $s$ with the maximum uncertainty. This data point is passed to the Generalized Benders Decomposition algorithm and the solution time (the label) is recorded, the labeled datapoint $\{s,y\}$ is appended in the training dataset, and the datapoint with label $s$ is removed form the pool. The surrogate model is trained again using the new training dataset and this loop continues until the maximum number of iterations is reached. The outcome of this approach is the surrogate model $\hat{m}$.

\begin{remark}
 \normalfont In this paper, we used a pool-based active learning approach with uncertainty-based sampling. However different active learning paradigms have been proposed. These paradigms can potentially be exploited for learning how to initialize decomposition-based and monolithic-based solution algorithms for different applications. For example, stream-based selective sampling can be used for learning online the surrogate models for the performance function, thus enabling the online learning of the optimal initialization. 
\end{remark}

\section{Application to mixed integer economic model predictive control for real-time operation of chemical processes} \label{application}

\begin{figure*}
    \centering
    \includegraphics[scale=0.5]{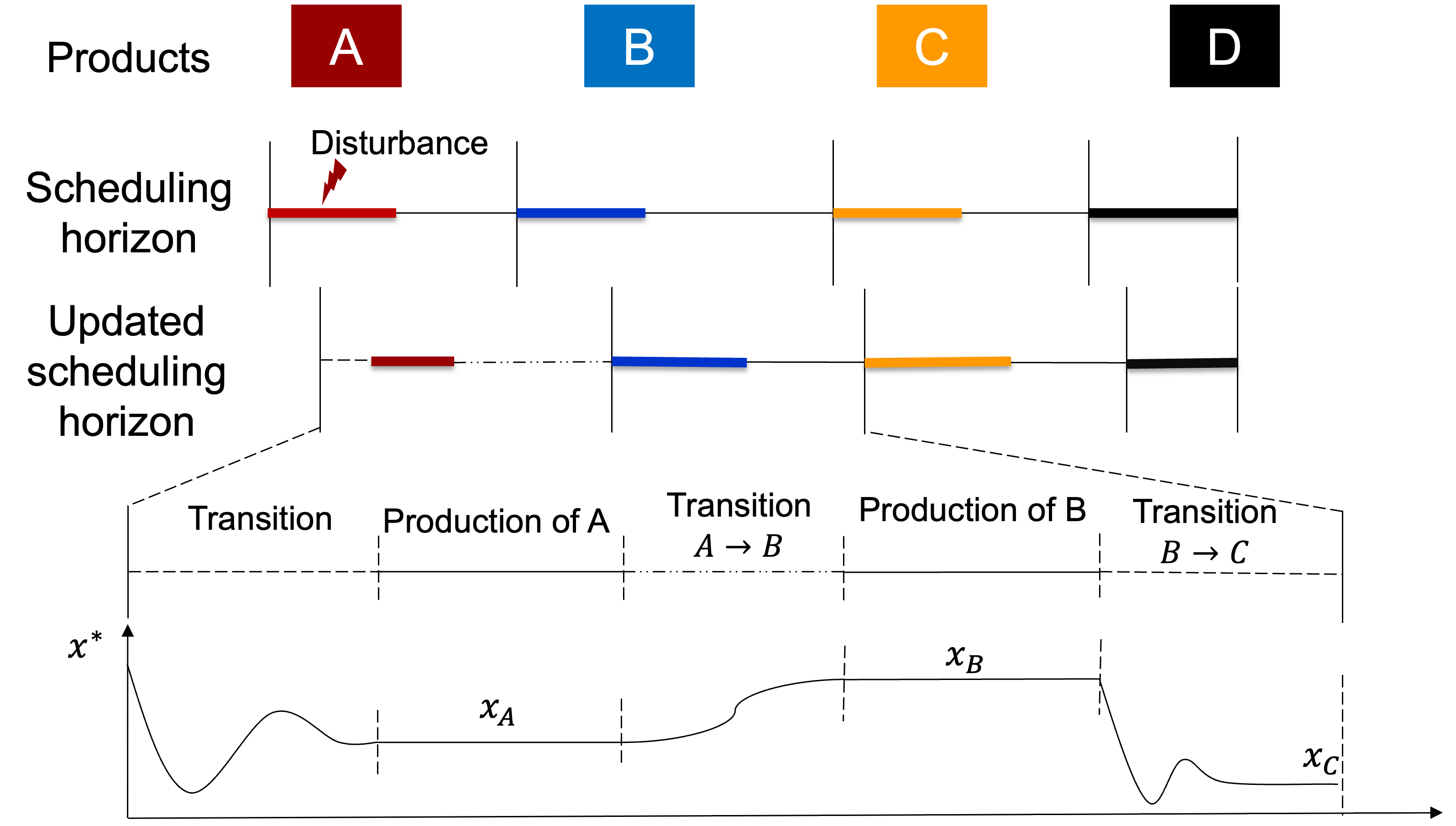}
    \caption{Schematic of rescheduling}
    \label{fig:MIP_MPC_schem}
\end{figure*}

In this section, we apply the proposed method for learning to initialize Generalized Benders Decomposition for the solution of mixed integer economic model predictive control problems that arise in the operation of chemical processes. Specifically, we consider a continuous manufacturing system whose operation can be affected by disturbances in the scheduling, e.g., change in demand, and control level, e.g., change in the inlet conditions of the process, as presented in Fig.~\ref{fig:MIP_MPC_schem}. Once a disturbance affects the system, an optimization problem is solved to determine the production sequence and the transitions between the products. In this section, we present the optimization model, the decomposition-based solution approach, the data generation process, and the evaluation of the different learning approaches.

\begin{algorithm}[t!]
\caption{Regression based initialization of Generalized Benders decomposition}\label{alg:online initialization}
\KwData{Surrogate model $\hat{m}$, Optimization problem $P$, value of parameters $p$}
\KwResult{Problem solution}
Compute the features of the problem $\nu(P(p))$\;
Determine the optimal number of cuts $n_{cuts}^* = \arg \min_{n \in \mathcal{N}_{cuts}} \hat{m}(N, \nu(P(p)))$\;
Add $n^*_{cuts}$ cuts in the master problem\;
Solve the optimization problem using Algorithm~1\;
\end{algorithm}

\begin{algorithm}[t!]
\caption{Learning surrogate model for solution time via active learning} \label{alg:AL algo}
\KwData{Number of evaluations $N$, initial labeled dataset $\mathcal{C} = \{s,y\}$, Pool of labels $\mathcal{C}_p$, Surrogate model $\hat{m}$}
\KwResult{Surrogate model $\hat{m}$}
Train surrogate model $\hat{m}$ on the initial dataset $\mathcal{C}$\;
\While{$i \leq N$}{
Select a data point with features $s$ from pool $\mathcal{C}_p$ based on maximum uncertainty sampling strategy $s \in \arg \max_{s \in \mathcal{C}_p} \sigma(s)$\;
Evaluate label $y$ for $s$ using GBD\;
Append data $\{s,y\}$ in set $\mathcal{C} = \mathcal{C} \cup \{s,y\}$\;
Remove data-point $s$ from pool $\mathcal{C}_p = \mathcal{C}_p \setminus \{s\}$\;
Train surrogate model $\hat{m}$ using set $\mathcal{C}$\;
$i=i+1$\;
}
\end{algorithm}

\subsection{Optimization model}
\subsubsection{Scheduling model}
We will consider the case where an isothermal CSTR is used to produce $N_p$ products over a time horizon of $H$ hours which is discretized into $N_s$ slots $(N_s=N_p)$. We will assume that while the system is following a nominal schedule a disturbance affects the system at time $T_0$, as presented in Fig.~\ref{fig:MIP_MPC_schem}. Under this setting, every slot $k$ except the first one has two regimes; a production regime where a product is manufactured and a transition regime where a transition occurs from the operating point of the product manufactured at slot $k$ to the operating point of the product manufactured at slot $k+1$. The first slot has three regimes; a transition regime that captures the transition from some intermediate state (where the system is due to the disturbance) to the operating point of the product manufactured in the first slot, a production regime, and another transition regime which considers the transition from the product manufactured at the first slot to the product manufactured in the second slot. 

We define a binary variable $W_{ik}$ which is equal to one if product $i$ is manufactured in slot $k$ and zero otherwise. Also, we define a binary variable $Z_{ijk}$ which is equal to one if a transition occurs from product $i$ to $j$ in slot $k$, and variable $\hat{Z}_i$ which is equal to one if a transition occurs from an intermediate state to product $i$ in the first slot. At every time point, only one product can be manufactured. The logic constraints are:
\begin{equation} \label{eq: logic cons sched}
    \begin{aligned}
        &\sum_{i =1}^{N_p} W_{ik} =1 \ \forall k=1,..,N_s\\
        & Z_{ijk} \geq W_{ik} + W_{j,k+1}-1 \ \forall i,j,k \neq {N_s} \\
        & \hat{Z}_i = W_{i1} \ \forall i
    \end{aligned}
\end{equation}
The starting time of slot $k$ is $T_k^s$, the ending time $T_k^e$, the production time of product $i$ in slot $k$ is $\Theta_{ik}$, and the transition time in slot $k$ is $\theta^t_{k}$. The timing constraints are:
\begin{equation}\label{eq: timing cons sched}
    \begin{aligned}
        T_k^e & = T_k^s + \sum_{i=1}^{N_p} \Theta_{ik} + \theta^t_{k} \ \forall k=1,.., N_s\\
        T_{k+1}^s &= T_k^e \ \forall k=1,..,N_s-1\\
        T_{N_s}^e &= H-T_0 \\
        \Theta_{ik} &\leq H W_{ik} \ \forall i=1,..,N_p, k=1,..,N_s\\
        \theta_{1}^t & = \sum_{i=1}^{N_p} \sum_{j=1}^{N_p} Z_{ij1} \theta_{ij1} + \sum_{i=1}^{N_p}\hat{Z}_i \hat{\theta}_i \\
        \theta_{k}^t &= \sum_{i=1}^{N_p} \sum_{j=1}^{N_p} Z_{ijk}\theta_{ijk} \ \forall k=2,..,N_s\\
        \theta_{ijk} &\geq \theta_{ij}^{min} \ \forall i,j,k\\
        \hat{\theta}_i &\geq \hat{\theta}^{min}_i \ \forall i
    \end{aligned}
\end{equation}
where $\theta_{ijk}$ is the transition time from product $i$ to $j$ in slot $k$, $\hat{\theta}_i$ is the transition time from an intermediate state to the steady state of product $i$, $\theta_{ij}^{min}$ is the minimum transition time from product $i$ to $j$, and $\hat{\theta}_{i}^{min}$ is the minimum transition time from the intermediate state to the steady state of product $i$. The production rate of product $i$ is $r_i$, the production amount of product $i$ in slot $k$ is $q_{ik}$, and the inventory of product $i$ in slot $k$ is $I_{ik}$. The production constraints are
\begin{equation} \label{eq: inv cons sched}
    \begin{aligned}
        I_{ik} & = I_{ik-1} + r_i \Theta_{ik} - S_{ik} \ \forall i=1,..,N_p, k=2,..,N_p\\
        I_{ik} & = I_{i}^{0} + r_i \Theta_{ik} - S_{ik} \ \forall i=1,..,N_p, k=1,
    \end{aligned}
\end{equation}
where $I^{0}_{i}$ is the initial inventory of product $i$. The demand of product $i$ is $d_i$ and the due date for every product is in the end of the time horizon. The demand satisfaction constraints are
\begin{equation} \label{eq: dem sat cons sched}
    S_{iN_s} \geq d_i \ \forall i=1,..,N_p
\end{equation}

\subsubsection{Dynamic model}
We will assume that the system is described by a system of differential equations
\begin{equation}
    \dot{x} = F(x,u),
\end{equation}
where $x \in \mathbb{R}^{n_x}$, $u \in \mathbb{R}^{n_u}$ are the state and manipulated variables, and $F: \mathbb{R}^{n_x} \times \mathbb{R}^{n_u} \mapsto \mathbb{R}^{n_x}$ are vector functions. We will consider simultaneously all transitions between the products and discretize the differential equations using the method of orthogonal collocation on finite elements (using $N_f$ finite elements and $N_c$ collocation points). We define state variable $x_{ijkfc}$ and manipulated variable $u_{ijkfc}$ for a transition from product $i$ to $j$ in slot $k$, finite element $f$ and collocation point $c$. We also define variable $\hat{x}_{ifc}$ and manipulated variable $\hat{u}_{ifc}$ for a transition from the intermediate state to product $i$ in finite element $f$ and collocation point $c$ in the first slot. The discretized differential equations for transitions between products are
\begin{equation} \label{eq: discr odes prod prod}
    \begin{aligned}
        x_{ijkfc} & = F_d (x_{ijkfc}, u_{ijkfc}, \theta_{ijk}) \ \forall i,j,k,f,c \\
        x_{ijk11} & = x_{i}^{ss} \ \forall i,j,k\\
        x_{ijkN_f N_c} &= x_{j}^{ss} \ \forall i,j,k\\
        u_{ijk11} &= u_{i}^{ss} \ \forall i,j,k\\
        u_{ijkN_f N_c} &= u_{j}^{ss} \ \forall i,j,k
    \end{aligned}
\end{equation}
and the equations for the transition from the intermediate state to product $i$ are:
\begin{equation} \label{eq: discr odes interm prod}
    \begin{aligned}
        \hat{x}_{ifc} & = F_d (\hat{x}_{ifc}, \hat{u}_{ifc}, \hat{\theta}_{i}) \ \forall i,l \\
        \hat{x}_{i11} & = x^{*} \ \forall i\\
        \hat{x}_{iN_f N_c} &= x_{i}^{ss} \ \forall i\\
        \hat{u}_{i11} &= u^{*} \ \forall i\\
        \hat{u}_{iN_f N_c} &= u_{i}^{ss} \ \forall i
    \end{aligned}
\end{equation}
where $x_{i}^{ss}$, $u_{i}^{ss}$ are the steady state values of the state and manipulated variables of product $i$. A detailed expression for the discretized differential equations can be found at \cite{mitrai2022multicut}.
\subsubsection{Objective function}
The objective function has three terms; the first is the profit $\Phi_1$, the second is the transition cost between products $\Phi_2$, and the third is the transition cost from the intermediate state $\Phi_3$. These terms are equal to:

\begin{equation}
    \begin{aligned}
        \Phi_1 & = \sum_{i,k} p_{ik} S_{ik} - C_{ik}^{oper} q_{ik} - C^{inv} I_{ik} - \sum_{ijk} C_{ijk}^{trans} Z_{ijk}\\
        \Phi_2 &= \sum_{ijk} Z_{ijk} \alpha_u \bigg( \sum_{fc} N_{f}^{-1} t_{ijfck}^d \Lambda_{c N_c} (u_{ijfck}-u_{j}^{ss})^2 \bigg) \\ 
        \Phi_3 &= \sum_{i} \hat{Z}_{i} \alpha_u \bigg( \sum_{fc} N_{f}^{-1} \hat{t}_{ifc}^d \Lambda_{c N_c} (\hat{u}_{ifc}-u_j^{ss})^2 \bigg) 
    \end{aligned}
\end{equation}
where $p_i,C_i^{op}$ are the price and operating cost of product $i$, $C^{inv}$ is the inventory cost, $C_{ij}^{tr}$ is the transition cost from product $i$ to $j$, $a_u$ is a weight coefficient, and $\Lambda$ is the collocation matrix.

\subsubsection{Mixed Integer Model Predictive Control problem and decomposition-based solution}
The mixed integer MPC problem is:
\begin{equation} \label{eq: MIP_MPC}
    \begin{aligned}
        P(p) : \minimize \ \ & \Phi_1 -\Phi_2 - \Phi_3 \\
        \subt \ \ & \text{Eq.}~\ref{eq: logic cons sched}, \ref{eq: timing cons sched}, \ref{eq: inv cons sched}, \ref{eq: dem sat cons sched}, \ref{eq: discr odes prod prod}, \ref{eq: discr odes interm prod}.
    \end{aligned}
\end{equation}
where $p = \{ \{d_{i}\}_{i=1}^{N_p}$, $\{I_i^0\}_{i=1}^{N_p}$,$T_0$, $\{\theta_{ij}^{min}\}_{i=1,j=1}^{N_p N_p}$,$ \{\hat{\theta}_i\}_{i=1}^{N_p}$, $\{r_i\}_{i=1}^{N_p}$, $\{x_i^{ss}\}_{i=1}^{N_p}$, $\{u_i\}_{i=1}^{N_p}, x^* \}$. This problem has three sets of parameters; the ones related to the scheduling part of the problem $\hat{p} =\{ \{d_{i}\}_{i=1}^{N_p}$, $\{I_i^0\}_{i=1}^{N_p}$, $T_0$, $\{\theta_{ij}^{min}\}_{i=1,j=1}^{N_p N_p}$,$ \{\hat{\theta}_i\}_{i=1}^{N_p}$, $\{r_i\}_{i=1}^{N_p} \}$, ones related to the dynamic behavior of the system for transitions between products $ \check{p}= \{ \{x_i^{ss}\}_{i=1}^{N_p}$, $\{u_i^{ss}\}_{i=1}^{N_p} \}$, and the ones related to the transition from the intermediate state to the steady state of the different products $\tilde{p} = x^*$.

We will rewrite the mixed integer economic MPC problem (Eq.~\ref{eq: MIP_MPC}) as follows
\begin{equation}
    \begin{aligned}
    \maximize \ \ & \ \ \Phi_1 (w) - \sum_{ijk} Z_{ijk} f_{dyn}^{ijk} (\omega_{ijk},\theta_{ijk}) - \sum_{i} \hat{Z}_i f_{dyn}^{i}(\hat{\omega}_{i},\hat{\theta}_{i}) \\
    \subt \ \ & \ \ g_{sched}(w, \theta_{ijk}, \hat{\theta}_i; \hat{p}) \leq 0 \\
    & \ \ g^{dyn}_{ijk}(\theta_{ijk}, \omega_{ijk}; \check{p}) \leq 0 \ \forall i,j,k \\
    & \ \ \hat{g}^{dyn}_i(\hat{\theta}_{i}, \hat{\omega}_{ijk}; \hat{p}) \leq 0 \ \forall i
    \end{aligned}
\end{equation}
where $w=\{ W_{ik}, Z_{ijk},\hat{Z}_i, T_k^s, T_k^, \Theta_{ik}, \theta_{k}^t, S_{ik} \}$ are scheduling variables, $\omega_{ijk} = \{ x_{ijkfc},u_{ijkfc} \}$ are variables associated with the dynamic behavior of the system for a transition from product $i$ to product $j$ in slot $k$, and $\hat{\omega}_{i}=\{\hat{x}_{ifc},\hat{u}_{ifc} \}$ are variables associated with the dynamic behavior of the system for a transition from the intermediate state to the product $i$. $g_{sched}$ are the scheduling constraints (Eq.~\ref{eq: logic cons sched}, \ref{eq: timing cons sched}, \ref{eq: inv cons sched}, \ref{eq: dem sat cons sched}), $g_{ijk}^{dyn}$ are the discretized differential equations for transitions between products (Eq.~\ref{eq: discr odes prod prod}), and $\hat{g}_i^{dyn}$ are the discretized differential equations for transition from the intermediate state to product $i$ (Eq.~\ref{eq: discr odes interm prod})

If we fix the scheduling variables $w$ and the transition times $\theta_{ijk},\hat{\theta}_i$ then the dynamic optimization problems for all the transitions can be solved independently. We define as $\phi_{ijk}$ the value function of the dynamic optimization problem for a transition from product $i$ to product $j$ in slot $k$. The dynamic optimization problem for this transition is 
\begin{equation}
    \begin{aligned}
    \phi_{ijk} (\theta_{ijk}): =\minimize_{\omega_{ijk}, \tilde{\theta}_{ijk}}  \ \ & \ \ f_{dyn}^{ijk} (\omega_{ijk},\tilde{\theta}_{ijk})  \\
    \subt \ \ & \ \ g_{dyn}(\tilde{\theta}_{ijk}, \omega_{ijk}) \leq 0 \\
    & \ \ \tilde{\theta}_{ijk} = \theta_{ijk} \ \ : \lambda_{ijk}.
    \end{aligned}
\end{equation}
where $\lambda_{ijk}$ is the Lagrangean multiplier for the equality constraint. Similarly, we define the value function for a transition from the intermediate state to product $i$, $\hat{\phi}_i$, and the dynamic optimization problem for this transition is
\begin{equation}
    \begin{aligned}
    \hat{\phi}_{i} (\hat{\theta}_{i}, \check{p}): =\minimize_{\hat{\omega}_i, \check{\theta}_i} \ \ & \ \ f_{dyn}^{ijk} (\hat{\omega}_{i},\check{\theta}_{i})  \\
    \subt \ \ & \ \ g_{dyn}(\check{\theta}_{i}, \hat{\omega}_{i}; \check{p}) \leq 0 \\
    & \ \ \check{\theta}_{i} = \hat{\theta}_{i} \ \ : \hat{\lambda}_{i}.
    \end{aligned}
\end{equation}

We note that these dynamic optimization problems are always feasible since the transition times are bounded from below by the minimum transition times, i.e., $\theta_{ijk} \geq \theta_{ij}^{min}$, $\hat{\theta}_i \geq \hat{\theta}_i^{min}$. The value functions $\phi_{ijk}, \hat{\phi}_i$ can be approximated with Benders cuts given by the following equations \cite{geoffrion1972generalized}
\begin{equation}
    \begin{aligned}
    \eta_{ijk} \geq & \phi_{ijk}(\bar{\theta}_{ijk}^l) - \lambda_{ijk}^l (\theta_{ijk} - \bar{\theta}_{ijk}^l)  \\
    \hat{\eta}_i \geq & \hat{\phi}_i (\bar{\hat{\theta}}_i^l) - \hat{\lambda}_i^l (\hat{\theta}_i - \bar{\hat{\theta}}_i^l),
    \end{aligned}
\end{equation}
where $l \in \mathcal{L}$ denotes the number of points used to approximate the value functions. The original problem can be reformulated as 
\begin{equation}
    \begin{aligned}
    \maximize \ \ & \ \ \Phi_1 (w) - \sum_{ijk} Z_{ijk} \eta_{ijk} - \sum_{i} \hat{Z}_i \hat{\eta}_{i} \\
    \subt \ \ & \ \ g_{sched}(w, \theta_{ijk}, \hat{\theta}_i) \leq 0 \\
    & \ \     \eta_{ijk} \geq \phi_{ijk}^l(\bar{\theta}_{ijk}^l) - \lambda_{ijk}^l (\theta_{ijk} - \bar{\theta}_{ijk}^l) \ \forall i,j,k,l  \\
    & \ \ \hat{\eta}_i \geq  \hat{\phi}_i (\bar{\hat{\theta}}_i) - \hat{\lambda}_i^l (\hat{\theta}_i - \bar{\hat{\theta}}_i^l) \ \forall i,l.
    \end{aligned}
\end{equation}
To solve the problem we use a hybrid multicut Generalized Benders Decomposition algorithm proposed in \cite{mitrai2022multicut}. In this algorithm the solution of the master problem provides the production sequence and the transition times, and Benders cuts are added to approximate the transition cost. In this case, the dynamic optimization problems between the products depend only on the transition time, whereas the transition from the intermediate state depends on the transition time and the concentration of the intermediate state. Therefore, the initialization of the algorithm considers only the optimal number of cuts added to approximate the transitions between the products, i.e. how many points should be used to approximate $\phi_{ijk}$. We use the same number of cuts for all transitions.

\begin{algorithm}[t!]
\caption{ Unlabeled data generation procedure for creating the pool of unlabeled data for active learning }\label{alg:create pool}
\KwData{ Number of data points $N_{data}$, Maximum number of discretization points $N_{cuts}$, upper and lower bounds for complicating variables $\theta_{ijk} \in [\theta_{ij}^{min}, 5\theta_{ij}^{min}$], Demand distribution for every product, Distribution of inlet concentration in the reactor, Scheduling horizon $H$ }
\KwResult{ Pool of unlabeled data $\mathcal{C}_p$}
$\mathcal{C}_p = \{ \}$\;
$i=0$\;
\While{$i\leq N_{data}$}{
Select at random a time point $T_0 \sim U(0,H)$\;
Introduce a disturbance in the inlet concentration of the rector $c_0 \sim U(0.8,1.2)$\;
Generate new demand for every product $d_i \sim U(\underbar{d}_i, \bar{d}_i)$\;
Compute inventory $I_0$ at time $T_0$\;
Get the value of the concentration in the reactor $x^*$\;
Obtain minimum transition times $\hat{\theta}_i^{min}$ from $x^*$ to $x_i^{ss}$\;
Form and solve the master problem\;
\eIf{Master problem is feasible}{
\For{$n=2:N_{cuts}$}{
Add $n$ cuts to add in the master problem\;
Compute features $s$ and append the data point in the pool $\mathcal{C}_p$\;}
}{The demand can not be satisfied at the end of the time horizon due to the disturbance, data point is not considered}
}
\end{algorithm}

\subsection{Application of active learning approach}
For the application of active learning, first, we generate the pool of unlabeled data as presented in Algorithm~\ref{alg:create pool}. We assume that at a random time point $T_0$ between 0 and the end of the scheduling horizon, the demand of all the products and the inlet concentration of the reactor change simultaneously based on some probability distributions. The statistics of the demand are presented in Table~\ref{demand distr} and we assume that the inlet concentration follows a uniform distribution with a low value of 0.8 and a high value of 1.2. Once the disturbance occurs, we compute the inventory, based on the realization of the initial schedule that was followed for $T_0$ hours and the concentration $x^*$ inside the reactor. Next, the minimum transition time from the intermediate state $x^*$ to the steady state $x_i^{ss}$ is computed for all the products. Given this information, we formulate a new master problem where the scheduling horizon is $H-T_0$ and check its feasibility. If the master problem is infeasible, then the disturbance that was generated can not be rejected such that the system can meet the demand at the end of the scheduling horizon. If the master problem is feasible, then different numbers of cuts are added to the master problem, as presented in lines 12-14 in Algorithm~\ref{alg:create pool}, and the features of the problem and the number of cuts are added in the unlabeled pool $\mathcal{C}_p$. Note that the feasibility of the mixed integer MPC problem can be determined by the feasibility of the master problem without any cuts since the subproblems are always feasible and the only source of infeasibility in the mixed integer MPC problem is due to the inability to satisfy the demand. 

The features of the problem, in this case, are the time point $T_0$, the concentration in the reactor $x^*$, the inlet flowrate $Q$, the demand of the products $\{d_i\}_{i=1}^{N_{prod}}$, inventory of the products $\{I_i^0\}_{i=1}^{N_{prod}}$, state of the system, i.e., production or transition, and the number of cuts added $n_{cuts}$. For the state of the system, we use one-hot encoding, i.e., $state \in \{0,1\}$. Overall the features $s_i$ for data point $i$ are:
\begin{equation*}
    s_i = [T_0,x^*,Q, \{d_i\}_{i=1}^{N_{prod}},\{I_i^0\}_{i=1}^{N_{prod}}, state, n_{cuts}]
\end{equation*}
These features form the pool $\mathcal{C}_p = \{ s_i \}_{i=1}^{N_{pool}}$ of data points with $N_{pool} = 49000$. Next, we generate a small number of data points $(N_{initial} = 10)$ and evaluate the CPU time. We allow 100 evaluations $N=100$, i.e., 100 data points are labeled. The computational time for obtaining the data is 726 seconds. The active learning appraoch is implemented using scikit-learn \cite{scikit-learn} and the mulitcut Generalized Benders Decomposition is implemented in Pyomo \cite{hart2017pyomo}, the master problem is solved using Gurobi \cite{gurobi} and the subproblem is solved with IPOPT \cite{wachter2006implementation}.

\subsection{Active learning results }
We compare the proposed active learning approach, denoted as GP-AL, with a random sampling of 110 points from the pool using different surrogates models such as Gaussian Process (GP) with Matern kernel, Neural Network (NN), Random Forest (RF) and Decision Tree (DT). The hyperparameters of the Gaussian Process, Random Forest, and the Decision Tree were set equal to their default values while the neural network has 3 layers, 150 neurons per layer, and \verb|tanh| as activation function. 

We generate 100 new disturbances that are not part of the pool and were not used for training in either the active or supervised learning approach. The solution time statistics for the different models are presented in Table~\ref{compare AL computational results}. From the results, we observe that the active learning approach leads on average to $66.5 \%$ reduction in CPU time compared to the standard application of Generalized Benders Decomposition, whereas the Gaussian process, Neural Network, Random Forest, and Decision Tree lead to $53\%, 43\%, 51\%$ and $33\%$ reduction respectively. 

Additionally, for the 100 random disturbances considered, the solution time obtained from the active learning approach is always lower than the solution time without the addition of cuts. The minimum percentage reduction in solution time for the active learning approach is $0.09 \%$, whereas for the surrogate models learned via random sampling, the minimum reduction is negative, i.e., the solution time of the proposed approach is higher than the original implementation of Generalized Benders Decomposition. These results indicate that the optimal number of cuts identified by optimizing the surrogate model trained via random sampling is highly suboptimal. Although these results can be justified by the fact that only 110 data points were used for training, they also highlight the importance of selecting the proper data points to label in cases where obtaining the labels is computationally expensive.

\subsection{Application of supervised learning}
We also consider the case where the labeling cost is not significant. Specifically, for every data point in the pool $\mathcal{C}_p$ with features $s_i$ generated in the previous section, we solve the optimization problem and obtain the solution time $y_i$, leading to a data set $\mathcal{D} = \{ s_i,y_i\}_{i=1}^{49000}$. We use this data set for training three surrogate modes; a Decision Tree, a Random Forest, and a Neural Network using scikit-learn \cite{scikit-learn}. For the Decision Tree and the Random Forest, we used the default values of the hyperparameters. The Neural Network had 3 layers with 150 neurons, the activation function was $\tanh$, the learning rate was equal to $10^{-4}$, and the regularization parameter $\alpha$ was set equal to 0.01. Note that in this case, we do not train Gaussian Processes since the dataset has 49000 data points and Gaussian Processes have cubic complexity on the number of samples.

\begin{figure}[t]
    \centering
    \resizebox{1.0\columnwidth}{!}{\includegraphics{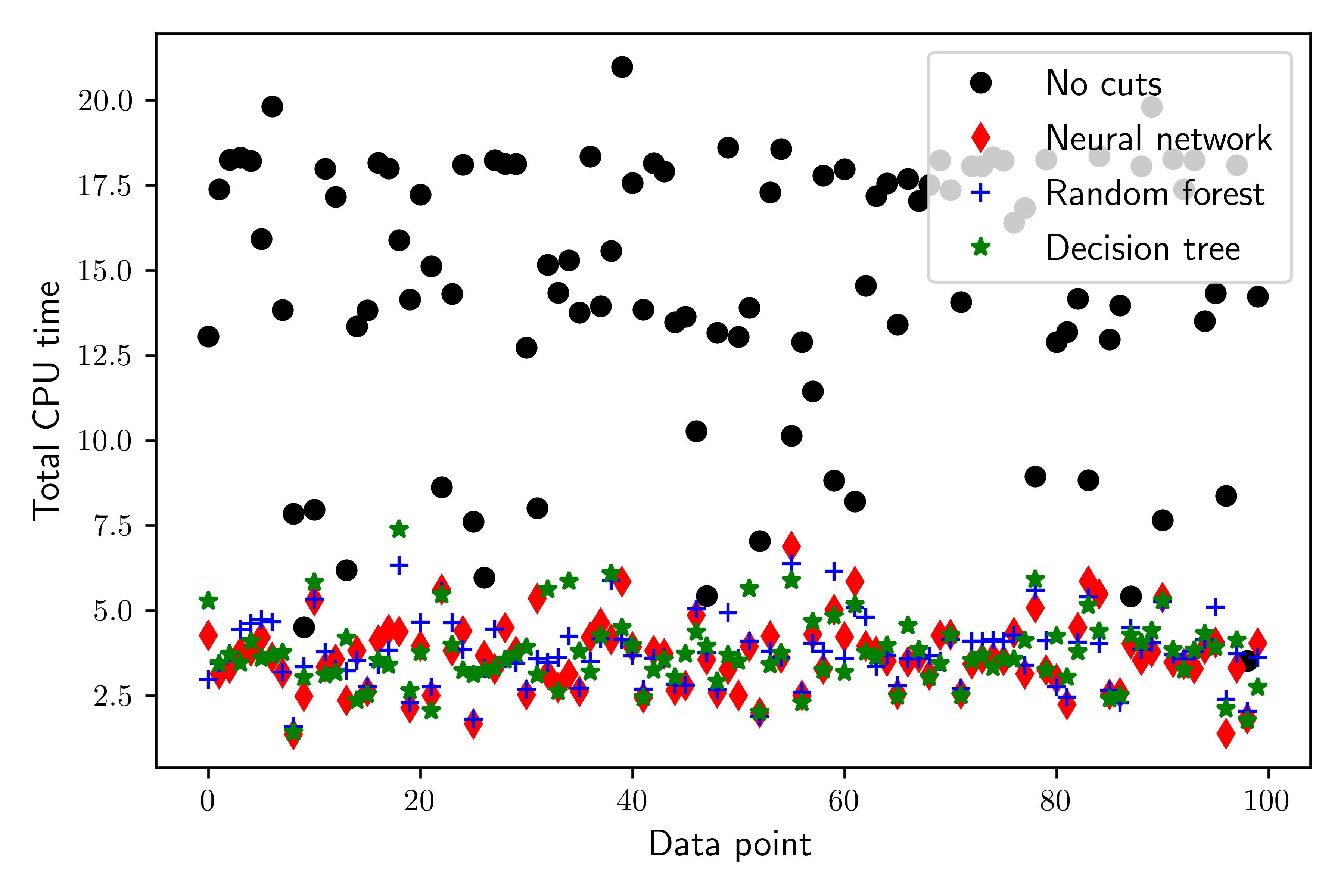}}
    \caption{Solution time of the proposed approach with different surrogate models for 49000 training data points.}
    \label{fig:sol_time_fig}
\end{figure}

Once the surrogate models were trained, we generate 100 random disturbances that change simultaneously the demands and the inlet concentration as in the previous section. These disturbances are different than the disturbances generated in the previous section. We solve the optimization problem for every disturbance by initializing the Generalized Benders Decomposition algorithm using the number of cuts suggested by optimizing the different surrogate models. The total CPU time for the different disturbances is presented in Fig.~\ref{fig:sol_time_fig} and the solution time statistics in Table~\ref{regression computational results}. From the results, we observe that the average total CPU time without the addition of cuts (No cuts) is $14.7$ seconds. The selection of the optimal number of cuts to add to the master problem leads to  a $70 \%$ reduction in CPU time. From the three surrogate models, the Neural Network shows the maximum improvement in total CPU time, although the average reduction is similar for all surrogate models. Furthermore, the minimum reduction is positive for all models, indicating that the solution time with the proposed initialization is lower than the original implementation of the algorithm without cuts. Finally, the time presented in Table~\ref{regression computational results} is the total time required to determine the optimal number of cuts and solve the problem. For case study considered, the time to determine the optimal number of cuts is in the order of $10^{-2}$ seconds for the decision tree and the neural network and in the order of $10^{-1}$ seconds for the random forest. Thus, the process of determining the number of cuts and adding them is significantly smaller than the solution time.

\begin{table}[t]
\centering
\caption{Computational time for the proposed approach for different surrogate models. NC refers to solving the problem without the addition of cuts in the first iteration.}
\resizebox{\columnwidth}{!}{
\begin{tabular}{ccccccc}
\hline
\multirow{2}{*}{\begin{tabular}[c]{@{}c@{}}Solution \\ statistics\end{tabular}} & \multicolumn{6}{c}{Initalization strategy}  \\ \cline{2-7} 
                 & NC      & GP-AL &  GP    & NN     & RF     & DT\\ \hline
Aver. CPU time   & 13.7    & 4.31  & 6.22   & 7.67   & 6.34   & 9.11 \\
Aver. red.       & -       & 66.5  & 53.44  & 43.52  & 51.61  & 33.64 \\
Aver. fold red.  & -       & 3.33  & 2.49   & 2.18   & 2.38   & 1.75     \\
Max. red. $(\%)$       & -       & 81.3  & 81.41  & 79.57  & 81.96  & 77.55    \\
Min. red. $(\%)$       & -       & 0.09  & -2.66  & -0.02  & -26.47 & -18.82\\
Max fold red.   & -        & 5.35  & 5.38   & 4.48   & 5.54   & 4.45 \\
Min fold red.    & -       & 1.00  & 0.97   & 0.99   & 0.79   & 0.84    \\ \hline
\end{tabular}}
\label{compare AL computational results}
\end{table}

\begin{table}[t]
\centering
\caption{Distribution of the demand}
\begin{tabular}{cccc}
\hline
\multirow{2}{*}{Product} & \multirow{2}{*}{Nominal value} & \multicolumn{2}{c}{Distribution (Uniform)} \\
                        &                                & low                  & high                \\ \hline
1                       & 600                            & -100                 & 100                 \\
2                       & 550                            & -15                  & 15                  \\
3                       & 600                            & -30                  & 30                  \\
4                       & 1200                           & -20                  & 20                  \\
5                       & 2000                           & -400                 & 400                 \\ \hline
\end{tabular}
\label{demand distr}
\end{table}

\begin{table}
\centering
\caption{Computational time for the proposed approach for different surrogate models trained via supervised learning}
\resizebox{\columnwidth}{!}{
\begin{tabular}{ccccc}
\hline
\multirow{2}{*}{\begin{tabular}[c]{@{}c@{}}Solution \\ statistics\end{tabular}} & \multicolumn{4}{c}{Initalization strategy}  \\ \cline{2-5} 
      & No cuts & \begin{tabular}[c]{@{}c@{}}Neural \\ networks\end{tabular} & \begin{tabular}[c]{@{}c@{}}Random \\ forests\end{tabular} & \begin{tabular}[c]{@{}c@{}}Decision\\  trees\end{tabular} \\ \hline
Average CPU time (sec)   & 14.7    & 3.63   & 3.78   & 3.74    \\
Average reduction $(\%)$   & -       & 71.71 & 70.50 & 70.53     \\
Average fold reduction   & -       & 4.23    & 4.01    & 4.10     \\
Max. reduction $(\%)$     & -       & 84.85   & 83.88     & 86.40     \\
Min. reduction $(\%)$  & -       & 25.66    & 17.30   & 21.13    \\
Max fold reduction  & -       & 6.60   & 6.20     & 7.35   \\
Min fold reduction    & -       & 1.34     & 1.20  & 1.26    \\ \hline
\end{tabular}}
\label{regression computational results}
\end{table}

\section{Conclusions and discusion} \label{concls}

The repeated solution of large-scale decision-making problems arises frequently in the operation of process systems. The efficient solution of such problems with monolithic optimization algorithms can be challenging, especially in online settings. Although decomposition-based solution methods have been widely used to solve large-scale optimization problems, their off-the-shelf implementation is nontrivial. In this work, we proposed a ML-based approach for optimally initializing cutting plane-based decomposition-based solution algorithms, such as Generalized Benders. We use active learning to guide the generation of labeled data for learning a surrogate model that predicts the solution time of Generalized Benders Decomposition for a given problem and the initial set of cuts added in the master problem. The proposed approach is applied to a case study of mixed integer economic model predictive control. The numerical results show that the optimal initialization of the algorithm can significantly reduce the solution time up to $70 \ \%$ and the computation of the optimal number of cuts using the learn surrogate model does not incur additional computational cost. These results highlight the ability of active learning to guide the data generation process for cases where obtaining the solution time is computationally expensive. 

\begin{figure*}[t!]
	\centering
	\includegraphics[scale=0.53]{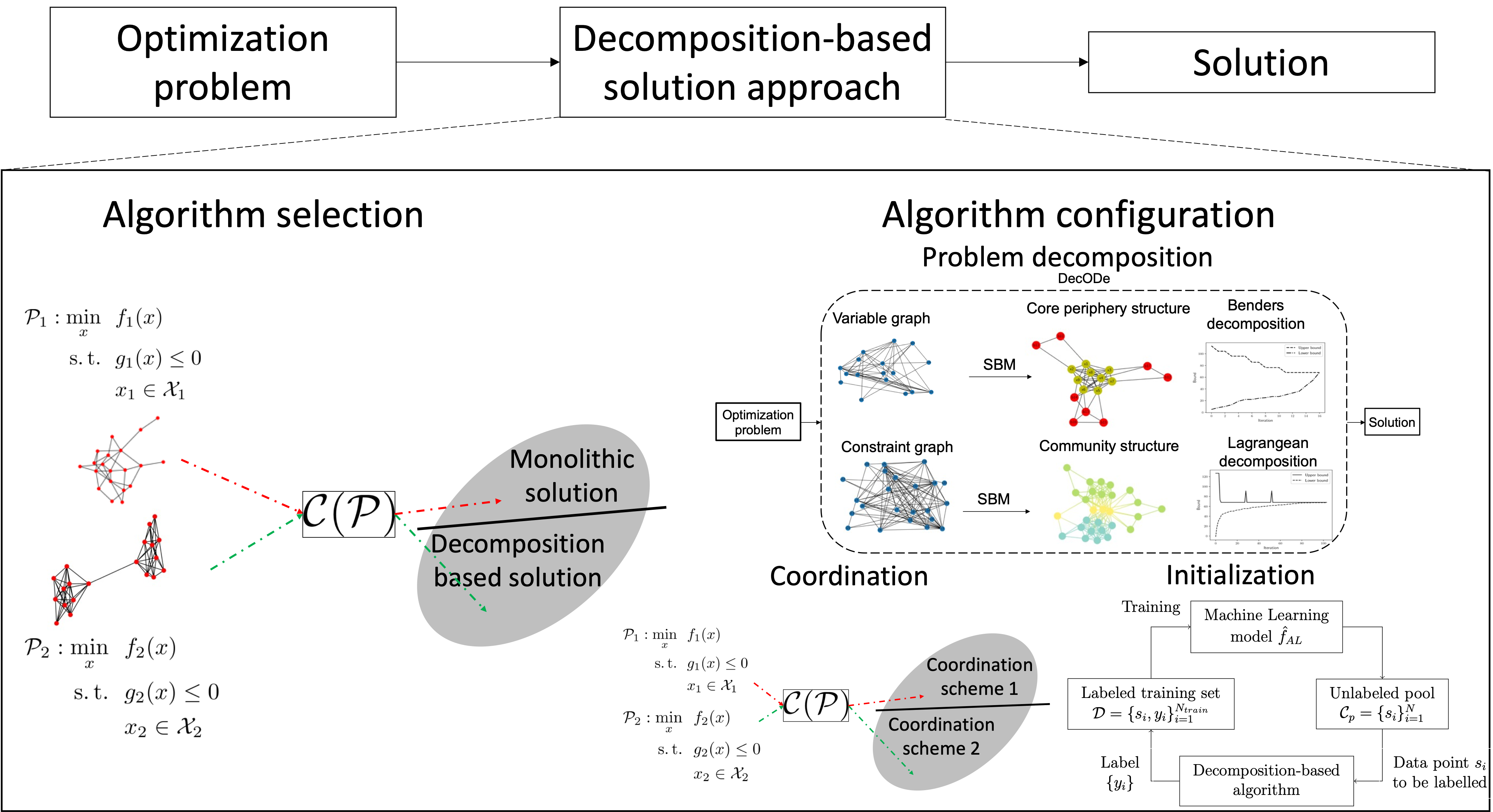}
	\caption{Framework for automated decomposition-based solution algorithm selection and configuration via artificial intelligence and network science}
    \label{fig: AI4OPt fram}
\end{figure*}

In general, the solution of an optimization problem with decomposition-based solution methods poses four questions: (1) Whether a decomposition-based method should be selected over a monolithic one, (2) how to decompose the optimization problem, (3) which coordination scheme should be used, and (4) how to initialize the algorithm.  The results presented in this paper and in \cite{mitrai2023whentodec}, in tandem with our previous work on learning the underlying structure of an optimization problem \cite{allman2019decode,mitrai2022stochastic, mitrai2021iecr}, show that artificial intelligence tools in combination with network science can be used to create an automated framework for decomposition-based solution algorithm selection and configuration.

The overall framework is presented in Fig.~\ref{fig: AI4OPt fram} and has two parts. The first focuses on algorithm selection, where given an optimization problem graph classification techniques can be used to determine whether, and potentially which, decomposition-based solution algorithm should be used. In the second part, the configuration of the selected decomposition-based solution algorithm is considered. The configuration has three aspects, problem decomposition, coordination, and initialization. For the problem decomposition, structure detection algorithms can be implemented to learn the underlying structure of the problem and use it as the basis for the application of the decomposition-based solution algorithm using DecODe \cite{mitrai2022stochastic}. 

Regarding the coordination, although it was not considered in this paper or in \cite{mitrai2023whentodec}, it is a critical aspect of decomposition-based solution. For hierarchical decomposition-based methods, the coordination is done via cuts and the different coordination schemes correspond to different cut generation and management techniques. For distributed algorithms, the coordination considers   different methods to update the dual variables or the Lagrangean multipliers. Consider for example the case of Lagrangean decomposition, where the the dual variables can be updated using subgradient, cutting plane, bundle, and trust region methods \cite{conejo2006decomposition}. The choice of which coordination scheme is the best, especially for mixed integer nonlinear programming problems is not obvious. ML tools such as classification and regression, possibly coupled with geometric deep learning, might aid the selection of the coordination scheme as presented in Fig.~\ref{fig: AI4OPt fram}. Finally, for the initialization, supervised learning can be used to learn surrogate models for the computational performance of the algorithm which can subsequently be used to optimally initialize the algorithm. 

This framework, and any ML-based approach, requires data for learning the parameters of the surrogate models used in the different tasks. Although multiple libraries of optimization problems exist \cite{minlplib,minlplibANDnlp,gleixner2021miplib}, these are unlabeled data since for a given problem the solution time for a given algorithm and a given configuration is not known. Building a library where not only the optimization problem but also the solution time and possibly other information regarding the solution is stored would significantly reduce the required time for building such ML-based approaches for decomposition-based solution algorithms. Finally, the computational effort required to obtain the labels can be reduced either via active learning or semi-supervised learning \cite{balestriero2023cookbook}.

\section{Acknowledgements}
Financial support from NSF-CBET is gratefully acknowledged. 

\bibliographystyle{elsarticle-num}
\bibliography{sample}

\end{document}